\documentclass[11pt]{amsart}

\usepackage[margin=1.5in]{geometry}

\usepackage{amssymb,amscd,amsxtra}
\usepackage{latexsym}
\usepackage[dvips]{graphics,epsfig}

\usepackage{graphicx}
\usepackage{exscale}
\usepackage{psfrag}
\usepackage{stmaryrd}

\usepackage{float}

\usepackage{amsmath,amsfonts,amssymb}
\usepackage{amsthm}
\usepackage{enumerate}
\usepackage{bbm}

\usepackage[dvips]{color}

\newcommand{\nA}{\mathbb A}
\newcommand{\nB}{\mathbb B}

\newcommand{\nR}{\mathbb R}

\newcommand{\cC}{\mathcal C}
\newcommand{\cQ}{\mathcal Q}
\newcommand{\cT}{\mathcal T}

\newcommand{\seg}{\mathcal S}

\newcommand{\la}{\langle}\newcommand{\ra}{\rangle}

\newcommand{\spn}{{\phi^2}}

\newtheorem*{theorem}{Theorem}

\theoremstyle{definition}
\newtheorem*{problem}{Problem}

\long\def\symbolfootnote[#1]#2{\begingroup\def\thefootnote{\fnsymbol{footnote}}
\footnote[#1]{#2}\endgroup}

\begin{document}

\begin{titlepage}

\rule{0pt}{0pt}

\vspace*{2cm}

{\Huge  \bf \sc  Ceva's triangle inequalities}

\vspace*{2cm}

{\LARGE \sc \'Arp\'ad B\'enyi\footnotemark[1]}

\vspace*{2mm}

{\normalsize Department of Mathematics, Western Washington
University,

516 High Street, Bellingham, Washington 98225, USA

E-mail: {\tt Arpad.Benyi@wwu.edu}}

\vspace*{1cm}

{\LARGE \sc Branko \'{C}urgus
}

\vspace*{2mm}

{\normalsize Department of Mathematics, Western Washington
University,

516 High Street, Bellingham, Washington 98225, USA

E-mail: {\tt Branko.Curgus@wwu.edu}}

\footnotetext[1]{Corresponding author. This work is partially supported by a grant from the Simons Foundation (No.~246024 to \'Arp\'ad B\'enyi).
}

\symbolfootnote[0]{{\bf 2010 Mathematics Subject Classification:} Primary 26D07, 51M04; Secondary 51M15, 15B48}

 \end{titlepage}

\title[Ceva's triangle inequalities]{}

\begin{abstract}
We characterize triples of cevians which form a triangle independent of the triangle where they are constructed. This problem is equivalent to solving a three-parameter family of inequalities which we call Ceva's triangle inequalities. Our main result provides the parametrization of the solution set.
\end{abstract}

\date{\today}


 \maketitle

\section{Introduction}

In this article we investigate a class of inequalities of the form
\begin{equation} \label{eqCineq}
\big| f(a,b,c,\rho) - f(b,c,a,\sigma) \bigr| < f(c,a,b,\tau) < f(a,b,c,\rho) + f(b,c,a,\sigma)
\end{equation}
where $(a, b, c)$ belongs to
\[
\cT = \bigl\{(a, b, c) \in {\mathbb R}^3 :
  | a - b\ \! | < c < a + b   \bigr\}.
\]
and $f:\cT \times \nR \to \nR_+$ is some function with positive values.
We will refer to \eqref{eqCineq} as {\em Ceva's triangle inequalities}.
Clearly, \eqref{eqCineq} expresses the fact that for some parameters $\rho, \sigma, \tau \in \nR$ the lengths
\[
f(a,b,c,\rho), \quad f(b,c,a,\sigma), \quad f(c,a,b,\tau)
\]
form a triangle whenever the lengths $a, b, c$ form a triangle. This problem, albeit much more general, is in the spirit of Bottema, Djordjevi\'c, Jani\'c, Mitrinovi\'c and Vasi\'c \cite[Chapter~13]{Bottema}, in particular the statement 13.3 which is attributed to Brown~\cite{Brown}.

For example, consider the function
\[
h(a,b,c,\rho) = \frac{\rho}{a} \sqrt{2(a^2b^2+b^2c^2+c^2a^2)-(a^4+b^4+c^4)}.
\]
Recall that $h(a,b,c,2)$ is the length of the altitude orthogonal to the side of length $a$ in a triangle with sides $a, b, c$. It is easy to see that, with $a=3, b=8, c=9$, the altitudes do not form a triangle. Thus, \eqref{eqCineq} fails to hold for all $(a,b,c) \in \cT$ with $\rho= \sigma = \tau$ for the function $h$. However, if we replace $f$ with $1/h$, then it turns out that~\eqref{eqCineq} does hold for all $\rho = \sigma = \tau \in \nR\setminus\{0\}$ and all $(a,b,c)\in \cT$. For a connection between the reciprocals of the altitudes and a Heron-type area formula, see Mitchell's note~\cite{Mitchell}.

In this article we will study \eqref{eqCineq} with
\begin{equation} \label{eqfunf}
f(a,b,c,\rho) = \sqrt{\rho(\rho-1)a^2 + \rho b^2+(1-\rho)c^2}.
\end{equation}
By Stewart's theorem, $f(a,b,c,\rho)$ gives the length of the {\em cevian}  $AA_{\rho}$ in a triangle $ABC$ with sides $BC=a,
C\!A=b$, $AB=c$ and where the point $A_\rho$ lies on the line $BC$ with $\overrightarrow{BA}_\rho =\rho \, \overrightarrow{BC}$. Similarly,
$f(b,c,a,\sigma) = BB_\sigma$ and $f(c,a,b,\tau) = CC_\tau$
where $B_\sigma$ and $C_\tau$ lie on the lines $C\!A$ and $AB$ respectively, with $\overrightarrow{CB}_\sigma=\sigma \overrightarrow{C\!A}$ and $\overrightarrow{AC}_\tau=\tau \, \overrightarrow{AB}$. Indeed, it is this connection with cevians that inspired the name of the inequalities~\eqref{eqCineq}. Now, it is well known that the medians of a triangle form a triangle as well. This implies that~\eqref{eqCineq} is satisfied for all $(a,b,c) \in \cT$ when $\rho = \sigma = \tau = 1/2$. In fact,~\eqref{eqCineq} is true for all $(a,b,c)\in \cT$ whenever $\rho = \sigma = \tau \in \nR$.  That is, for any triangle $ABC$ and for all $\xi\in\nR$, the cevians $AA_\xi, BB_\xi,$ and $CC_\xi$ always form a triangle. This statement is proved in  the book by Mitrinovi\'c, Pe\v{c}ari\'c and Volenec~\cite[Chapter~I.3.14]{Mitrinovic} where it is attributed to Klamkin~\cite{Klamkin}; see also the articles of Hajja~\cite{Hajja1, Hajja2} for a geometric proof.

\begin{figure}[p]

\psfrag{A}[][]{\begin{picture}(0,0)
            \put(-2,0){\makebox(0,0)[l]{$A=C_\tau$}}
                        \end{picture}}
\psfrag{B}[][]{\begin{picture}(0,0)
            \put(-14,-15){\makebox(0,0)[l]{$B$}}
                        \end{picture}}
\psfrag{C}[][]{\begin{picture}(0,0)
            \put(-5,1){\makebox(0,0)[l]{$C$}}
                        \end{picture}}

\psfrag{Cr}[][]{\begin{picture}(0,0)
            \put(-5,3){\makebox(0,0)[l]{}}
                        \end{picture}}

\psfrag{Bq}[][]{\begin{picture}(0,0)
            \put(-6,0){\makebox(0,0)[l]{$B_\sigma$}}
                        \end{picture}}

\psfrag{Ap}[][]{\begin{picture}(0,0)
            \put(-7,1){\makebox(0,0)[l]{$A_\rho$}}
                        \end{picture}}

\psfrag{App}[][]{\begin{picture}(0,0)
            \put(-9,0){\makebox(0,0)[l]{$A_\rho^\prime$}}
                        \end{picture}}

\setlength{\abovecaptionskip}{-2pt}%
\setlength{\belowcaptionskip}{0pt}%

\resizebox{0.7\linewidth}{!}{%
  \includegraphics{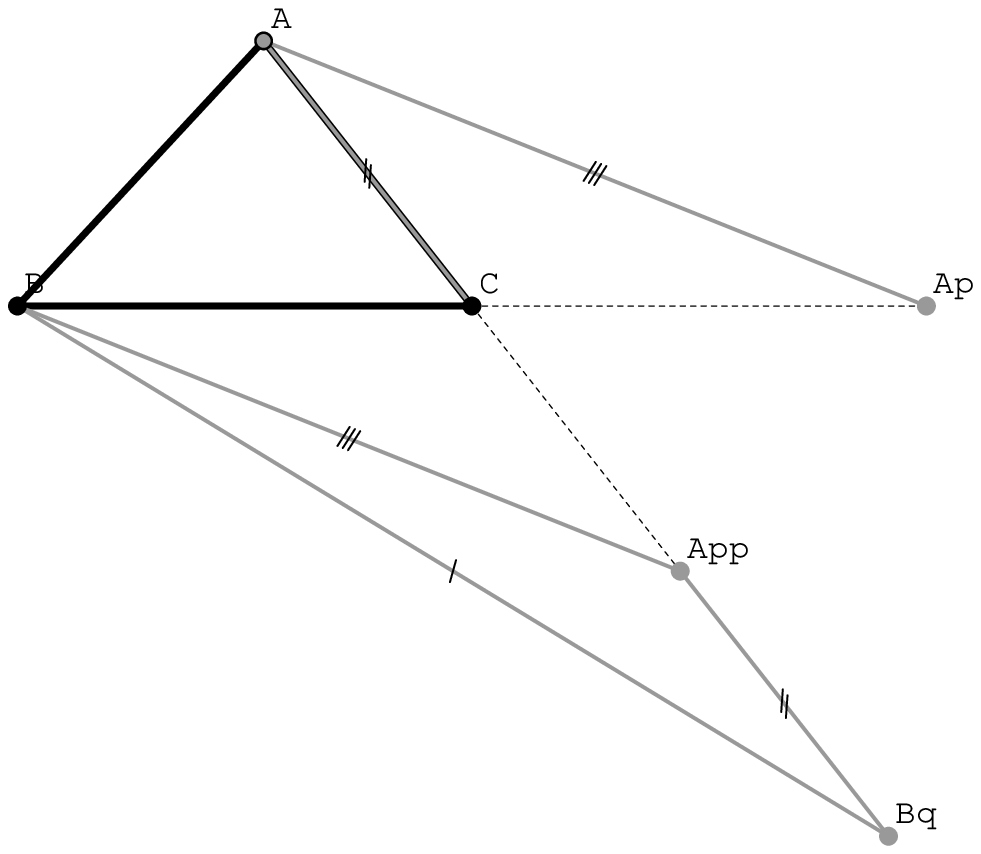}}
    \caption{Cevians always form a triangle}
\label{uniex}

\end{figure}

Interestingly, for an arbitrary triangle $ABC$, the triple $(\rho, \sigma, \tau)=(2,-2,0)$ produces the sides $AA_{2},
BB_{-2}$, and $CC_{0}$ which also form a triangle. In other words, for the given $f:\cT\times \nR\to \nR_+$, \eqref{eqCineq} also holds for all $(a,b,c)\in \cT$ and some non-diagonal triple $(\rho, \sigma, \tau)\not\in\bigl\{(r, s, t)\in \nR^3: r=s=t\bigr\}.$ This can be proved by using the methods followed in Figure~\ref{uniex}.

The property does not hold for all triples of real numbers. Consider, for example, $(\rho, \sigma, \tau)=(1/4, 1/2, 5/6)$. In the right triangle $ABC$ with $AB=1, BC=\sqrt{8}, C\!A=3$, we easily
find $AA_{\rho}=\sqrt{3/2}\approx 1.225$, $BB_{\sigma}=3/2=1.5$, and
$CC_\tau=17/6 \approx 2.833$. Thus, $AA_\rho+BB_\sigma<CC_\tau$,
which means that the three cevians do not form a triangle in
this case. However, for the same triple but the right triangle
$ABC$ with $AB=1, BC=1, C\!A=\sqrt{2}$, we find
$AA_\rho=\sqrt{17}/4\approx 1.031, BB_\sigma=\sqrt{2}/2\approx
0.707,$ and $CC_\tau=\sqrt{37}/6\approx 1.014$, which are now
clearly seen to form a triangle.

\setlength{\abovecaptionskip}{-2pt}%
\setlength{\belowcaptionskip}{0pt}%

\begin{figure}[p]

\psfrag{tau}[][]{\begin{picture}(0,0)
            \put(25,0){\makebox(0,0)[l]{\LARGE $\tau$}}
                        \end{picture}}
\psfrag{sigma}[][]{\begin{picture}(0,0)
            \put(0,25){\makebox(0,0)[l]{\LARGE $\sigma$}}
                        \end{picture}}
\psfrag{rho}[][]{\begin{picture}(0,0)
            \put(-35,0){\makebox(0,0)[l]{\LARGE $\rho$}}
                        \end{picture}}

\resizebox{0.7\linewidth}{!}{%
  \includegraphics{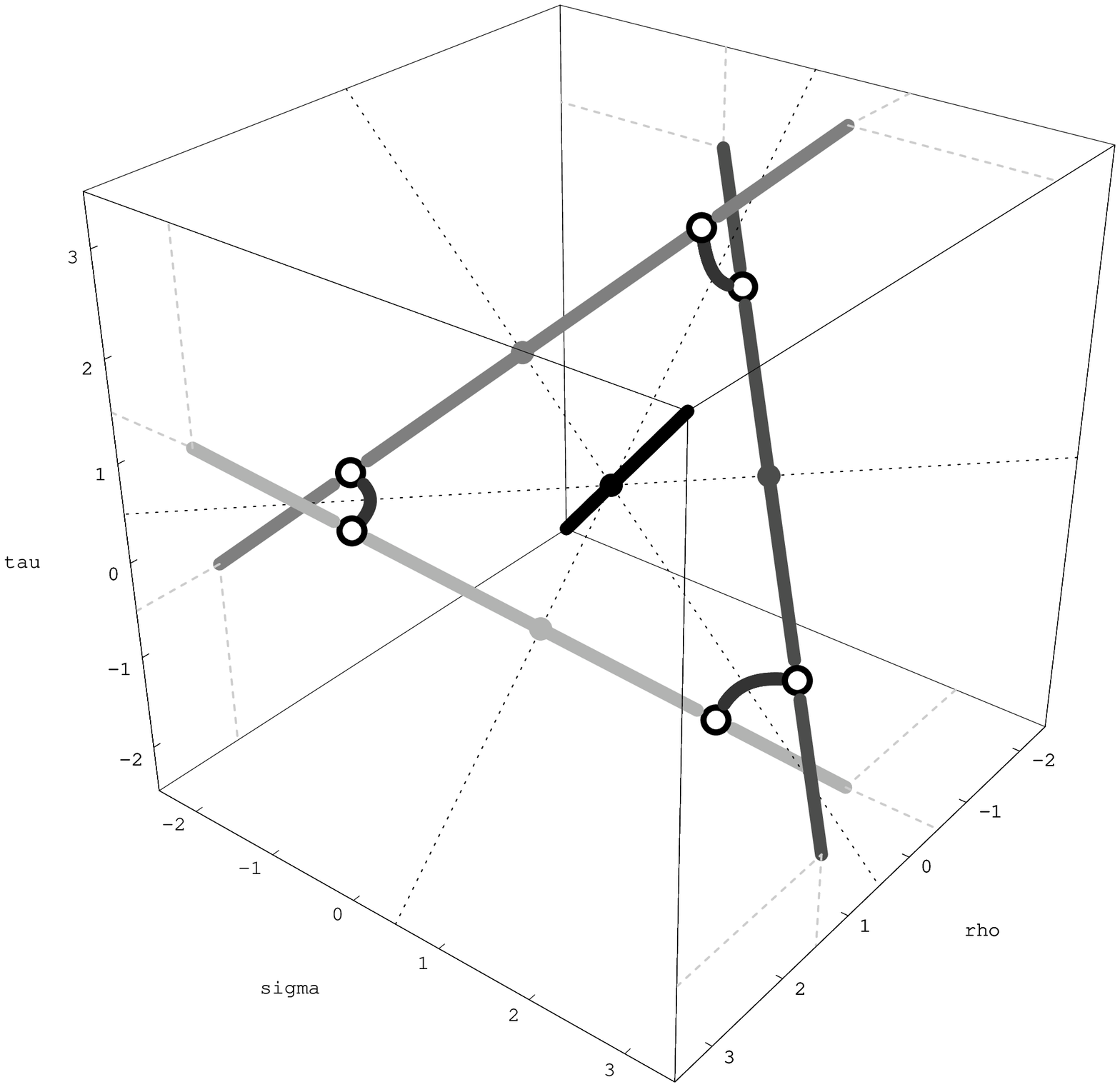}}
    \caption{The set $\nA$}
\label{SetA}

\end{figure}

The cautionary examples above lead us to the following problem.

\begin{problem}
Characterize the set $\mathbb A \subset \nR^3$ of all triples $(\rho, \sigma, \tau)\in {\mathbb R}^3$ such that, for the function $f$ given by \eqref{eqfunf}, Ceva's triangle inequalities \eqref{eqCineq} hold for all $(a, b, c)\in \cT$. In other words, characterize the set of all triples $(\rho, \sigma, \tau)$ such that for {\em all} non-degenerate  triangles \,\!\! $ABC$ the cevians $AA_\rho, BB_\sigma$ and $CC_\tau$ form a non-degenerate triangle as well.
\end{problem}

The main result of this paper is the complete parametrization of the set $\nA$.

\begin{theorem}
Let $\phi = (1+\sqrt{5})/2$ denote the golden ratio. The set $\nA$ is the union of the following three sets
\begin{gather} \label{eqsetD}
\Bigl\{ (\xi, \xi, \xi) \ : \ \xi \in \nR \Bigr\}, \\
 \label{eqsetE}
 \Bigl\{(-\xi, 2-\xi, \xi),
 \ (\xi, -\xi, 2-\xi), \ (2-\xi, \xi, -\xi) \ : \  \xi \in \nR\setminus\{-\phi^{-1},\phi\} \Bigr\} \\
\intertext{and}
  \label{eqsetF}
 \Bigl\{ \bigl(\tfrac{1}{1-\xi}, 1\! -\! \tfrac{1}{\xi}, \xi \bigr), \
   \bigl( \xi, \tfrac{1}{1-\xi}, 1\! -\!\tfrac{1}{\xi} \bigr), \
    \bigl( 1\! -\!\tfrac{1}{\xi},  \xi, \tfrac{1}{1-\xi} \bigr) \ :
     \ \xi \! \in \!
 \bigl( \phi, \phi^2 \bigr)
  \Bigr\}.
\end{gather}
\end{theorem}

A picture of the solution set $\nA$ is given in Figure~\ref{SetA}. Notice that the following six points are excluded in the second set:
\begin{equation} \label{eq6pts}
\begin{array}{lll}
    \bigl(\phi^{-1}, \phi^2, -\phi^{-1} \bigr), & \quad
    \bigl( -\phi^{-1}, \phi^{-1}, \phi^2 \bigr), & \quad
    \bigl(\phi^{2}, -\phi^{-1}, \phi^{-1} \bigr),  \\[6pt]
    \bigl(-\phi, \phi^{-2}, \phi \bigr), & \quad
     \bigl(\phi,  - \phi, \phi^{-2} \bigr), & \quad
    \bigl(\phi^{-2}, \phi, -\phi\bigr).
 \end{array}
\end{equation}
These points are the only common accumulation points of both the second and the third set in the parametrization of $\nA$. This claim follows from the identities
\begin{gather*}
\left(-\phi, 2 - \phi, \phi \right) = \left(\tfrac{1}{1-\phi},  1 - \tfrac{1}{\phi}, \phi \right), \\
\bigl(-(-\phi^{-1}), 2 -(-\phi^{-1}), -\phi^{-1} \bigr) = \left(1 - \tfrac{1}{\phi^2},  \phi^2, \tfrac{1}{1-\phi^2}\right),
\end{gather*}
which are consequences of $\phi^2 - \phi - 1 = 0$. As we will see in Subsection~\ref{sdeg}, these six points are exceptional since for an arbitrary triangle they correspond to degenerate triangles formed by the three cevians.

The remainder of the paper is devoted to the proof of the Theorem.
First, in Subsection~\ref{1}, we introduce a one-parameter family of inequalities whose solution set $\nB$ is guaranteed to contain $\nA$. In Section~\ref{secB} we use a combination of arguments from analytic and synthetic geometry to prove that $\nB$ is a subset of the union of the sets defined in \eqref{eqsetD},  \eqref{eqsetE}, \eqref{eqsetF} and \eqref{eq6pts}. In Section~\ref{secA}, we use linear algebra to prove that the union of the sets in \eqref{eqsetD},  \eqref{eqsetE}, \eqref{eqsetF} is in fact contained in $\nA$ and that the points in \eqref{eq6pts} are not in $\nA$. This provides the parametrization given in the Theorem.

It is worthwhile noting that the discussion about the set $\nA$ reduces to the observation that an uncountable family ${\mathcal M}$ of $3\times 3$ matrices has a common invariant cone; that is, for a particular cone $\cQ$ we have $M\cQ \subset \cQ$ for all $M \in {\mathcal M}$.  The question about when does a finite family of matrices share an invariant cone is of current interest and already non-trivial. For example, in the ``simplest'' case of a finite family that is simultaneously diagonalizable one has a characterization for the existence of such an invariant cone, but this criteria does not seem to be easily applicable; see Rodman, Seyalioglu, and Spitkovsky \cite[Theorem 12]{Rodman} and also Protasov \cite{Protasov}.

\section{The set $\nB$}\label{secB}

Recall that, for a given $(a,b,c) \in \cT$
and a triple $(\rho, \sigma, \tau)\in {\mathbb R}^3$, we have the following expressions for the lengths of the cevians:
 \[
AA_\rho=f(a, b, c, \rho),\quad BB_\sigma=f(b, c, a, \sigma), \quad CC_\tau=f(c, a, b, \tau).
 \]
Specifically,
\begin{align*} 
AA_\rho & = \sqrt{\rho(\rho-1)a^2 + \rho b^2+(1-\rho)c^2}, \\
BB_\sigma & = \sqrt{(1-\sigma)a^2 + \sigma(\sigma-1)b^2 + \sigma c^2}, \\
CC_\tau & = \sqrt{\tau a^2 + (1-\tau)b^2 + \tau(\tau - 1)c^2}.
\end{align*}
Clearly, a triple $(\rho, \sigma, \tau)$ belongs to $\nA$ if and only if
the following {\it three-parameter family of inequalities} is satisfied:
\begin{equation}\label{triangle-cev}
 |AA_\rho-BB_\sigma| < CC_\tau < AA_\rho+BB_\sigma \quad \text{for all} \quad  (a, b, c)\in \cT.
\end{equation}

\subsection{The limiting inequalities}\label{1}

We now introduce the {\it one-parameter family of inequalities} emerging from~\eqref{triangle-cev} by letting $(a,b,c)$ belong to the faces of the closure of the infinite tetrahedron $\cT$.

Let $a=b+c,\, b, c > 0$. Then
 \[
AA_\rho=\sqrt{\rho(\rho-1)a^2+\rho (a-c)^2+(1-\rho)c^2}=\sqrt{(\rho a-c)^2}=|\rho a-c|.
 \]
Similarly,
 \[
BB_\sigma=|\sigma b-a|, \quad CC_\tau=|\tau c+b|.
 \]
Thus, for all $b, c>0$, we have
 \[
\bigl||\rho (b+c)-c|-|\sigma b-b-c|\bigr| \leq
 |\tau c+b|\leq
  |\rho (b+c)-c|+|\sigma b-b-c|.
 \]
Equivalently, by letting $b/c=t$, we have
\begin{equation}\label{limit1}
\bigl| |\rho (t+1)-1|-|\sigma t-t-1| \bigr| \leq
 |\tau +t|\leq
  |\rho (t+1)-1|+|\sigma t-t-1|,
\end{equation}
for all $t>0$.

Analogously, by considering the other two faces of $\cT$, we obtain the following two inequalities
\begin{equation}\label{limit2}
\bigl||\rho +t|-|\sigma (t+1)-1|\bigr| \leq
 |\tau t-t-1| \leq
  |\rho +t|+|\sigma (t+1)-1|,
\end{equation}
\begin{equation}\label{limit3}
 \bigl||\rho t-t-1|-|\sigma +t|\bigr| \leq
  |\tau (t+1)-1|\leq
   |\rho t-t-1|+|\sigma +t|,
\end{equation}
for all $t>0$.

\subsection{The definition and a partition of the set $\nB$}

By $\nB$ we denote the set of all $(\rho,\sigma,\tau) \in \nR^3$ which satisfy all three limiting inequalities~\eqref{limit1},  \eqref{limit2}, and~\eqref{limit3} for all $t>0$. By $\nB_\tau$ we denote the horizontal section of $\nB$ at the level $\tau$, that is, the subset of $\nB$ with a fixed $\tau$. Then $\nB = \bigcup_{\tau \in \nR} \nB_\tau$.

In the previous subsection we have shown that $\nA \subseteq \nB$. In the remainder of this section we study the sets $\nB_\tau$.

The following four points play the central role in the characterization:
\begin{alignat*}{2}
P_{11}(\tau) &:= (\tau, \tau, \tau), \qquad  & P_{12}(\tau) &:= (-\tau, 2-\tau, \tau), \\
P_{21}(\tau) &:= (2-\tau, -2+\tau, \tau), \qquad & P_{22}(\tau) &:= (2+\tau, -\tau, \tau).
\end{alignat*}
Note that the points $P_{11}(\tau), \tau \in \nR$, appear in \eqref{eqsetD} and most of the points $P_{12}(\tau)$, $P_{21}(\tau)$, $P_{22}(\tau)$, $\tau \in \nR$, appear in \eqref{eqsetE}.

\subsection{$\nB_\tau$ is a subset of the union of two lines}\label{3}

Let $\tau \in \nR$ be fixed. Let $\ell_j(\tau), j \in \{1,2\}$, be the lines determined by the points $P_{j1}(\tau)$ and $P_{j2}(\tau)$.

We first show that if $(\rho, \sigma, \tau) \in \nB_\tau$, then
\begin{equation} \label{eqlines}
|\rho (-\tau+1)-1|=|-\sigma\tau+\tau-1|.
\end{equation}
We break our discussion into three cases, depending on the range of the parameter $\tau$. If $\tau\leq 0$, and we let $t\rightarrow (-\tau)^{+}$ in~\eqref{limit1}, then the leftmost inequality gives~\eqref{eqlines}.
If $\tau\in (0, 1]$, then we let $t\rightarrow \Bigl(\displaystyle\frac{1-\tau}{\tau}\Bigr)^{-}$ in the first part of~\eqref{limit3} to get
 \[
\Bigl|\rho\frac{1-\tau}{\tau} - \frac{1-\tau}{\tau}-1\Bigr| = \Bigl|\sigma+\frac{1-\tau}{\tau}\Bigr|,
 \]
which implies~\eqref{eqlines}.
Finally, if $\tau > 1$, let $t=\displaystyle\frac{1}{\tau-1}$ in the first  inequality of~\eqref{limit2} to obtain
 \[
\Bigl|\rho+\frac{1}{\tau-1}\Bigr|=\Bigl|\sigma \Bigl(\frac{1}{\tau-1}+1\Bigr)-1\Bigr|,
 \]
which implies~\eqref{eqlines}.

Elementary considerations yield that the graph of~\eqref{eqlines} is the union $\ell_1(\tau)\cup \ell_2(\tau)$. Thus
\[
\nB_\tau \subseteq \ell_1(\tau)\cup \ell_2(\tau).
\]

\subsection{The sets $\nB_0$ and $\nB_1$}\label{stau01}
We show now that $\nB_0 = \bigl\{P_{ij}(0) \,|\, i, j\in \{ 1, 2 \} \bigr\}$. Recall that the points $(\rho, \sigma, \tau)$ must satisfy~\eqref{eqlines}. Substituting $\tau=0$ in this equation gives
$|\rho-1|=1$, thus $\rho=0$ or $\rho=2$.  If $\rho=0$, then~\eqref{limit2} is
 \[
\bigl||t|-|\sigma (t+1)-1|\bigr| \leq  t+1 \leq
  t+|\sigma(t+1)-1|, \quad t>0.
 \]
If we let $t\rightarrow 0^{+}$, we get $|\sigma-1|\geq 1\geq |\sigma-1|$, thus $|\sigma-1|=1$, that is $\sigma=0$ or $\sigma=2$. So our points are $(0, 0, 0)$ and $(0, 2, 0)$. If $\rho=2$, then~\eqref{limit3}  is
 \[
\bigl||t-1|-|\sigma+t|\bigr| \leq  1 \leq |t-1|+|\sigma+t| , \quad t>0.
 \]
Substituting $t=1$ gives $|\sigma+1|\geq 1\geq |\sigma+1|$, that is $\sigma=0$ or $\sigma=-2$. This gives the points $(2, 0, 0)$ and $(2, -2, 0)$. All in all, the four points found are exactly $P_{ij}(0),i, j\in \{ 1, 2 \}$.

Similarly, it can be shown that $\nB_1 = \bigl\{P_{ij}(1)\,|\,  i, j\in \{ 1, 2 \} \bigr\}$.

\subsection{A family of crosses associated with the limiting inequalities} \label{5}

In the rest of this section we assume that $\tau \in \mathbb \nR \setminus \{0, 1\}$.  In this subsection we establish some preliminary facts which are then needed in the last subsection to prove that
\begin{equation*} 
\nB_\tau = \bigl\{P_{ij}(\tau)\,|\,  i, j\in \{ 1, 2 \} \bigr\}
 \ \ \text{or} \ \
 \nB_\tau = \bigl\{P_{ij}(\tau)\,|\,  i, j\in \{ 1, 2 \} \bigr\} \cup \bigl( \ell_1(\tau)\cap \ell_2(\tau)\bigr).
\end{equation*}

Let $t > 0$ be fixed. Notice that inequalities~\eqref{limit1}, \eqref{limit2} and~\eqref{limit3} can be written in the form
\begin{equation*} 
\bigl||a(t)\rho-b(t)|-|c(t)\sigma-d(t)|\bigr|
 \leq |g(t)| \leq
  |a(t)\rho-b(t)|+|c(t)\sigma -d(t)|,
\end{equation*}
where $a, b, c, d, g$ are all affine functions of $t$; one of the coefficients of $g$ depends on the fixed value $\tau$. Since $g(t) = 0$ yields~\eqref{eqlines}, we will only consider the case $g(t) \neq 0$.

\begin{figure}[H]

\setlength{\abovecaptionskip}{0pt}%
\setlength{\belowcaptionskip}{0pt}%

\noindent
\begin{minipage}{0.49\linewidth}
\resizebox{\linewidth}{!}{
  \includegraphics{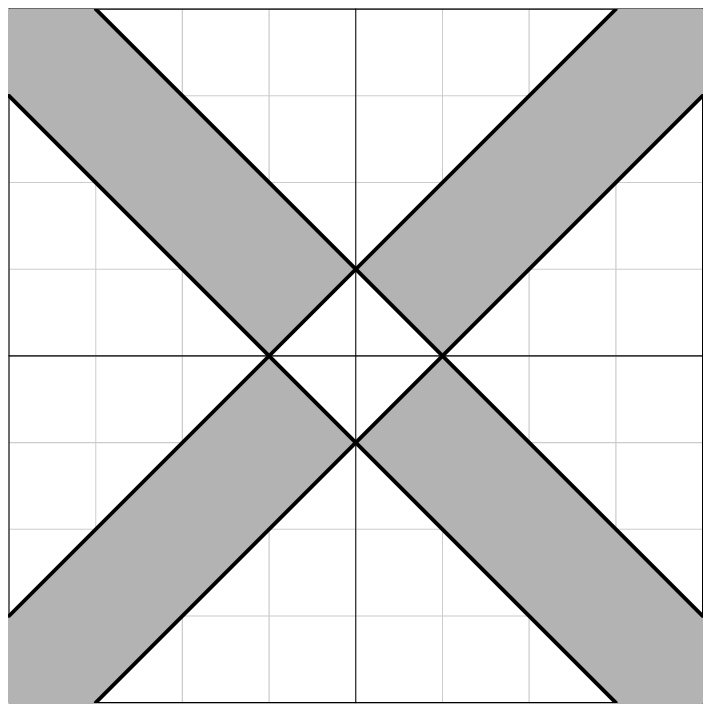}}
  \caption{The canonical cross} \label{piCross}
\end{minipage}
\begin{minipage}{0.49\linewidth}
\resizebox{\linewidth}{!}{
  \includegraphics{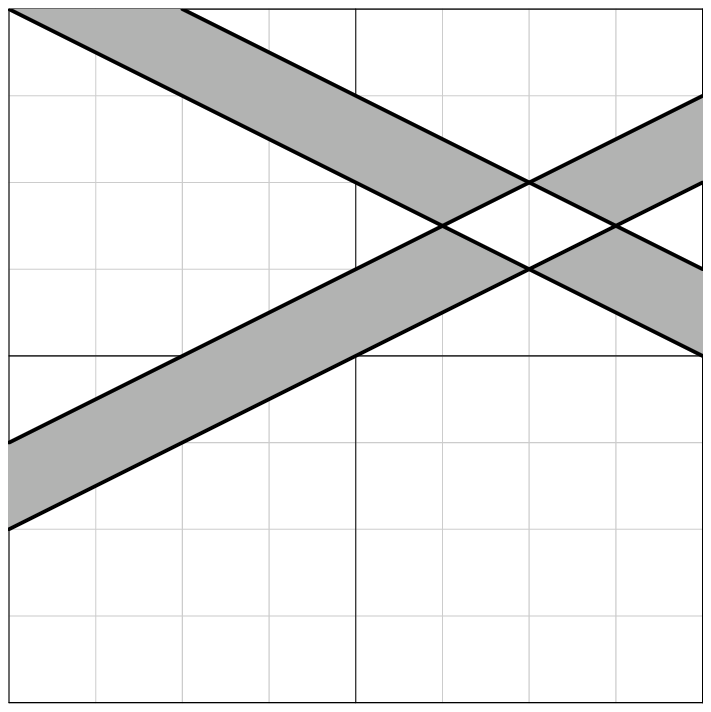}}
 \caption{Shifted and skewed} \label{piCrossS}
\end{minipage}

\end{figure}

Renaming $\rho=x$, $\sigma=y$, we are interested in fully understanding the geometric representation in the $xy$-plane of inequalities in the generic form
\begin{equation}\label{cross}
\bigl||ax-b|-|cy-d|\bigr| \leq  |g| \leq
 |ax-b|+|cy-d|,
\end{equation}
where $a, b, c, d, g \in\mathbb R\setminus \{0\}$. Note that~\eqref{cross} is equivalent to
\begin{equation}  \label{re-cross}
\biggl|\Bigl|\frac{a}{g}\Bigr|\Bigl|x - \frac{b}{a}\Bigr| - \Bigl|\frac{c}{g}\Bigr|\Bigl|y - \frac{d}{c}\Bigr|\biggr|
 \leq  1 \leq
 \Bigl| \frac{a}{g}\Bigr|\Bigl|x - \frac{b}{a}\Bigr| + \Bigl|\frac{c}{g}\Bigr|\Bigl|y-\frac{d}{c}\Bigr|.
\end{equation}
The canonical form of such inequalities is
\begin{equation}\label{standard-cross}
\bigl||x|-|y|\bigr| \leq  1\leq |x|+|y|.
\end{equation}
Elementary considerations show that the solution set of~\eqref{standard-cross} is the set represented in Figure~\ref{piCross}. We will  refer to it as the {\it canonical cross}. It is centered at the origin and its marking points are, in counterclockwise direction, $(1, 0), (0, 1), (-1, 0),$ and $(0, -1)$.  Note also that there are two pairs of parallel lines that define the boundary of the canonical cross, and the slopes of the lines that are not parallel have values $1$ and $-1$.

The change of coordinates
 \[
x \mapsto \bigl|\tfrac{a}{g}\bigr|\bigl(x-\tfrac{b}{a}\bigr), \qquad
 y \mapsto \bigl|\tfrac{c}{g}\bigr|\bigl(y-\tfrac{d}{c}\bigr),
 \]
transforms~\eqref{standard-cross} into \eqref{re-cross} and the canonical cross is transformed into a shifted and skewed version of it, illustrated  in Figure~\ref{piCrossS}. Hence, a shifted and skewed cross is the solution set of~\eqref{cross}.  Its center is at $\bigl(\frac{b}{a}, \frac{d}{c}\bigr)$, and its marking points are
$\bigl(\tfrac{b}{a}+\bigl|\tfrac{g}{a}\bigr|, \tfrac{d}{c}\bigr),$ 
 $\bigl(\tfrac{b}{a}, \tfrac{d}{c}+\bigl|\tfrac{g}{c}\bigr|\bigr),$ 
 $\bigl(\tfrac{b}{a}-\bigl|\tfrac{g}{a}\bigr|, \tfrac{d}{c}\bigr)$ 
  and 
 $\bigl(\tfrac{b}{a}, \tfrac{d}{c}-\bigl|\tfrac{g}{c}\bigr|\bigr).$
The interior of the parallelogram with vertices at the marking
points will be referred to as the {\it nucleus of the cross}.
The equations of two pairs of parallel boundary lines of this cross are
\begin{equation*}
y =  \tfrac{a}{c} \bigl(x - \tfrac{b}{a}\bigr) + \tfrac{d}{c} \pm \tfrac{g}{c} \quad \text{and} \quad y = - \tfrac{a}{c} \bigl(x - \tfrac{b}{a} \bigr) + \tfrac{d}{c} \pm \tfrac{g}{c}.
\end{equation*}

For a  given $t > 0$, the limiting inequality~\eqref{limit1} has the form~\eqref{cross}. Hence its solution set is a cross in the $\rho\sigma$-plane at the {\em fixed} level $\tau$. We denote this cross by $\mbox{Cross}_{\ref{limit1}}(t,\tau)$. The equations of the two pairs of the parallel boundary lines of this cross are
\begin{equation} \label{eqbls1}
y =  \tfrac{t+1}{t}  \bigl(x - \tfrac{1}{t+1} \bigr)\! + \tfrac{t+1}{t} \pm \tfrac{\tau+t}{t} \ \ \text{and} \ \ y = - \tfrac{t+1}{t}  \bigl(x - \tfrac{1}{t+1} \bigr)  + \tfrac{t+1}{t} \pm \tfrac{\tau+t}{t}.
\end{equation}
Similarly,  $\mbox{Cross}_{\ref{limit2}}(t,\tau)$ and $\mbox{Cross}_{\ref{limit3}}(t,\tau)$
refer to the solution sets of~\eqref{limit2} and~\eqref{limit3}.
The equations of the two pairs of the corresponding parallel boundary lines of these cross, respectively, are
\allowdisplaybreaks{\begin{alignat}{2} \label{eqbls2}
y & =  \tfrac{1}{t+1}  \bigl(x + \tfrac{t}{1} \bigr)\! + \tfrac{1}{t+1} \pm \tfrac{\tau t - t - 1}{t+1} \ \ & \text{and} \ \ y & = - \tfrac{1}{t+1}  \bigl(x + \tfrac{t}{1} \bigr)\! + \tfrac{1}{t+1} \pm \tfrac{\tau t - t - 1}{t+1}, \\ \label{eqbls3}
y & =  \tfrac{t}{1}  \bigl(x - \tfrac{t+1}{t} \bigr)\! + \tfrac{-t}{1} \pm \tfrac{\tau(t+1) - 1}{1} \ \ & \text{and} \ \ y & = - \tfrac{t}{1}  \bigl(x - \tfrac{t+1}{t} \bigr)\! + \tfrac{-t}{1} \pm \tfrac{\tau(t+1) - 1}{1}.
\end{alignat}
}
\!\!By $\mathcal C(\tau)$ we denote the family of all these crosses:
\[
{\mathcal C}(\tau) = \bigl\{  \mbox{Cross}_{\ref{limit1}}(t,\tau): \ t >0  \bigr\}
 \cup  \bigl\{  \mbox{Cross}_{\ref{limit2}}(t,\tau), \ t >0  \bigr\}
   \cup   \bigl\{  \mbox{Cross}_{\ref{limit3}}(t,\tau), \ t >0  \bigr\}.
\]

In the next subsection we will need the following two facts about ${\mathcal C}(\tau)$:

\begin{enumerate}[{\bf {\rm\bf  Fact} 1.}]
\item  \label{f1}
The points $P_{i1}(\tau), P_{i2}(\tau), i \in\{1,2\}$, lie on the distinct parallel boundary lines of {\em each} of the crosses in ${\mathcal C}(\tau)$.

\item  \label{f2}
For each slope $m \in \nR \setminus \{0,1\}$ there is a cross in $\cC(\tau)$ with the boundary lines going through $P_{11}(\tau)$ and $P_{12}(\tau)$ having slope $m$ and, consequently, with the boundary lines through $P_{21}(\tau)$ and $P_{22}(\tau)$ having slope $-m$.

\end{enumerate}

Fact\!{~\ref{f1} is verified by substitution of the $\rho\sigma$-coordinates of $P_{ij}(\tau), i,j \in\{1,2\}$, in the equations of the boundary lines of the crosses in ${\mathcal C}(\tau)$.
For example, the $\rho\sigma$-coordinates of $P_{11}(\tau)$ satisfy the left equation in \eqref{eqbls1} with $-$ sign. The table below gives the complete account:

\begin{center}
\begin{tabular}{|c|c|c|c|c|}\hline
 & $P_{11}(\tau)$ & $P_{12}(\tau)$ & $P_{21}(\tau)$ & $P_{22}(\tau)$ \\\hline
\eqref{eqbls1}
    & left eq. with $-$
      & left eq. with $+$
        & right eq. with $-$
          & right eq. with $+$ \\\hline
\eqref{eqbls2}
    & left eq. with $+$
      & left eq. with $-$
        & right eq. with $+$
          & right eq. with $-$  \\\hline
\eqref{eqbls3}
    & right eq. with $+$
       & right eq. with $-$
         & left eq. with $+$
            & left eq. with $-$ \\\hline
\end{tabular}
\end{center}

To prove Fact\!{~\ref{f2} we distinguish three cases. If $m>1$, then, by \eqref{eqbls1} and the table above,  $\mbox{Cross}_{\ref{limit1}}\bigl(1/(m-1),\tau\bigr)$ has the desired property. If $0 < m < 1$, then, by \eqref{eqbls2} and the table, $\mbox{Cross}_{\ref{limit2}}\bigl(-1+1/m,\tau\bigr)$ proves Fact\!{~\ref{f2}. If $m < 0$, then, by \eqref{eqbls3} and the table, $\mbox{Cross}_{\ref{limit3}}\bigl(-m,\tau\bigr)$ is the one whose existence is claimed in Fact\!{~\ref{f2}.

Notice that when the slope $m$ is that of the line $\ell_1$, the cross as  in Fact 2 is degenerate, that is, it is the union of the lines $\ell_1$ and $\ell_2$.

\subsection{Intersection of crosses determined by four points}\label{7}

By Fact~\ref{f1} we know that the four points $P_{ij}(\tau), i,j \in \{1,2\}$, belong to all the crosses in ${\mathcal C}(\tau)$. Moreover, there is possibly only one extra point having the property that it lies on all the crosses. That point is
\[
P_0(\tau) := \Bigl(\frac{1}{1-\tau}, 1-\frac{1}{\tau}, \tau \Bigr),
\]
which is the intersection of the lines $\ell_1(\tau)$ and $\ell_2(\tau)$ defined in Subsection~\ref{3}. Note that some of the points $P_0(\tau)$, $\tau \in \nR\setminus\{0,1\}$, appear in \eqref{eqsetF}.

Now, for which $\tau$-s could one expect $P_0(\tau)$ to belong to all the crosses in $\cC(\tau)$? Letting $t\rightarrow 0^{+}$ in~\eqref{limit1} and~\eqref{limit3}, we get
\[
\bigl||\rho-1|-1\bigr| \leq  |\tau| \leq
 |\rho-1|+1
\quad \text{and} \quad
\bigl||\sigma|-1\bigr| \leq  |\tau-1| \leq
  1+|\sigma|.
\]
If we substitute $\rho=1/(1-\tau)$ and $\sigma=(\tau-1)/\tau,$ we obtain that $\tau$ must satisfy
\[
\bigl||\tau|-|\tau-1|\bigr| \leq  |\tau^2-\tau| \leq
 |\tau|+|\tau-1|.
\]
Straightforward calculations show that the solution of the above inequalities is the following union of disjoint intervals:
\begin{equation*}
\Bigl[\tfrac{-1-\sqrt{5}}{2}, \tfrac{1-\sqrt{5}}{2} \Bigr] \cup \Bigl[\tfrac{3-\sqrt{5}}{2}, \tfrac{-1+\sqrt{5}}{2} \Bigr] \cup
\Bigl[\tfrac{1+\sqrt{5}}{2}, \tfrac{3+\sqrt{5}}{2} \Bigr],
\end{equation*}
or in terms of the golden ratio $\phi$:
\begin{equation*} 
  \bigl[ -\phi, -\phi^{-1} \bigr] \! \cup \!
  \bigl[ \phi^{-2}, \phi^{-1} \bigr]\! \cup\!
  \bigl[ \phi, \phi^2 \bigr].
\end{equation*}

In Subsection~\ref{sdeg} we will indeed show that, for $\tau$ in the interior of this set,  the point $P_0(\tau)$ belongs to all the crosses in the family $\cC(\tau)$.

Since $\tau \in \nR \setminus \{0,1\}$ is fixed, henceforth in this proof we will not emphasize the dependence on $\tau$. We denote by $\seg_j$ the open line segment determined by $P_{j1}$ and $P_{j2}, j \in \{1,2\}$.

We will break our discussion into three separate cases that take into account the relative positions of the five points $P_0, P_{11}, P_{12}, P_{21}, P_{22}$. We will repeatedly use Fact~\ref{f2} without explicitly citing it.

\begin{enumerate}[{\bf {\rm\bf  Case} 1.}]
\item  \label{c1}
If $\{P_0\} = \seg_1 \cap \seg_2$, then $\nB_\tau = \bigl\{P_{11}, P_{12}, P_{21}, P_{22} \bigr\}$.

\item  \label{c2}
If $P_0 \not\in \seg_1 \cup \seg_2$, then $\nB_\tau = \bigl\{P_{11}, P_{12}, P_{21}, P_{22} \bigr\}$.

\item  \label{c3}
If $P_0 \in \bigl(\seg_1 \setminus \seg_2\bigr) \cup \bigl(\seg_2 \setminus \seg_1\bigr)$, then $\nB_\tau = \bigl\{P_{11}, P_{12}, P_{21}, P_{22}, P_0 \bigr\}$.
\end{enumerate}

Notice that the line $\ell_1$ containing $P_{11}$ and $P_{12}$ has slope $m = 1-1/\tau \in \nR\setminus \{0,1\}$ and that the line $\ell_2$ containing $P_{21}$ and $P_{22}$ has slope $-m$. Note also that when $\tau = 1/2$ we have $m = -1$. In this case $P_0$ is the midpoint of $\seg_2$ and it is outside of $\seg_1$. That is, $m= -1$ can occur only in {Case}\!{~\ref{c3}.}

\smallskip

\noindent
{\bf Case}\!{\bf~\ref{c1}.} Any point $X$ that lies in $\seg_1\cup
\seg_2$ can be eliminated by a cross that has the boundary lines through $P_{11}$ and $P_{12}$ of slope $-m$, and boundary
lines through points $P_{21}$ and $P_{22}$ of slope $m$; see Figure~\ref{c1p1}. To eliminate a point $X$ outside of the segments $\seg_1$ and $\seg_2$, consider the cross having boundary lines that pass through $P_{11}$ and $P_{12}$ of slope $m-\epsilon$, and boundary lines
passing through $P_{21}$, $P_{22}$ of slope $-m+\epsilon$, with $\epsilon \in (0, m)$ sufficiently small; see Figure~\ref{c1p2}.

\begin{figure}[p]

\setlength{\abovecaptionskip}{-3pt}%
\setlength{\belowcaptionskip}{6pt}%

\noindent
\begin{minipage}{0.49\linewidth}
\psfrag{P11}[][]{\begin{picture}(0,0)
            \put(20,0){\makebox(0,0)[l]{$P_{11}$}}
                        \end{picture}}
\psfrag{P12}[][]{\begin{picture}(0,0)
            \put(20,-1){\makebox(0,0)[l]{$P_{12}$}}
                        \end{picture}}
\psfrag{P21}[][]{\begin{picture}(0,0)
            \put(-12,0){\makebox(0,0)[l]{$P_{21}$}}
                        \end{picture}}
\psfrag{P22}[][]{\begin{picture}(0,0)
            \put(20,0){\makebox(0,0)[l]{$P_{22}$}}
                        \end{picture}}

\psfrag{P0}[][]{\begin{picture}(0,0)
            \put(18,0){\makebox(0,0)[l]{$P_0$}}
                        \end{picture}}

\psfrag{X}[][]{\begin{picture}(0,0)
            \put(5,7){\makebox(0,0)[l]{$X$}}
                        \end{picture}}

\resizebox{\linewidth}{!}{%
  \includegraphics{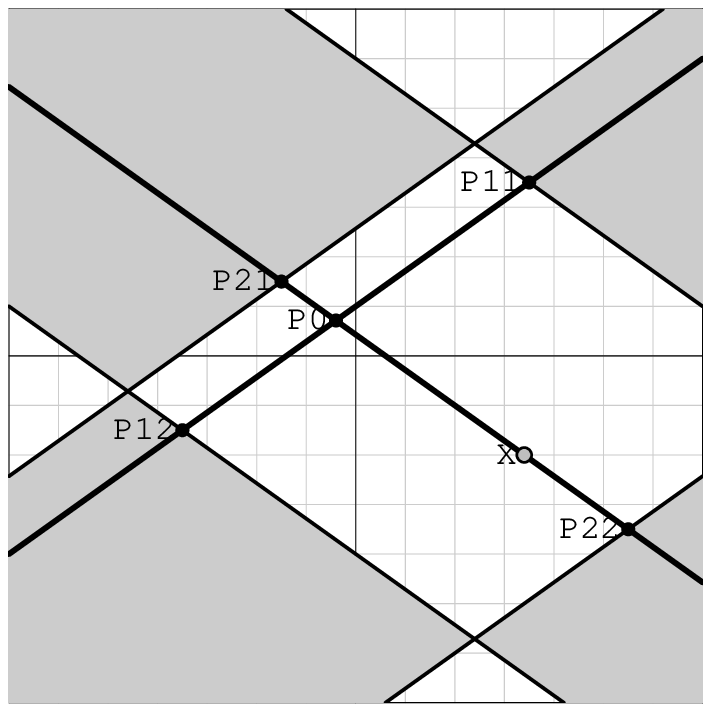}}
  \caption{} \label{c1p1}
\end{minipage}
\begin{minipage}{0.49\linewidth}

\psfrag{P11}[][]{\begin{picture}(0,0)
            \put(-3,8){\makebox(0,0)[l]{$P_{11}$}}
                        \end{picture}}
\psfrag{P12}[][]{\begin{picture}(0,0)
            \put(10,-8){\makebox(0,0)[l]{$P_{12}$}}
                        \end{picture}}
\psfrag{P21}[][]{\begin{picture}(0,0)
            \put(9,9){\makebox(0,0)[l]{$P_{21}$}}
                        \end{picture}}
\psfrag{P22}[][]{\begin{picture}(0,0)
            \put(-7,-4){\makebox(0,0)[l]{$P_{22}$}}
                        \end{picture}}

\psfrag{P0}[][]{\begin{picture}(0,0)
            \put(-4,-34){\makebox(0,0)[l]{$P_0$}}
            \put(2,-28){\vector(1, 4){6}}
                        \end{picture}}

\psfrag{X}[][]{\begin{picture}(0,0)
            \put(6,-6){\makebox(0,0)[l]{$X$}}
                        \end{picture}}

\resizebox{\linewidth}{!}{%
  \includegraphics{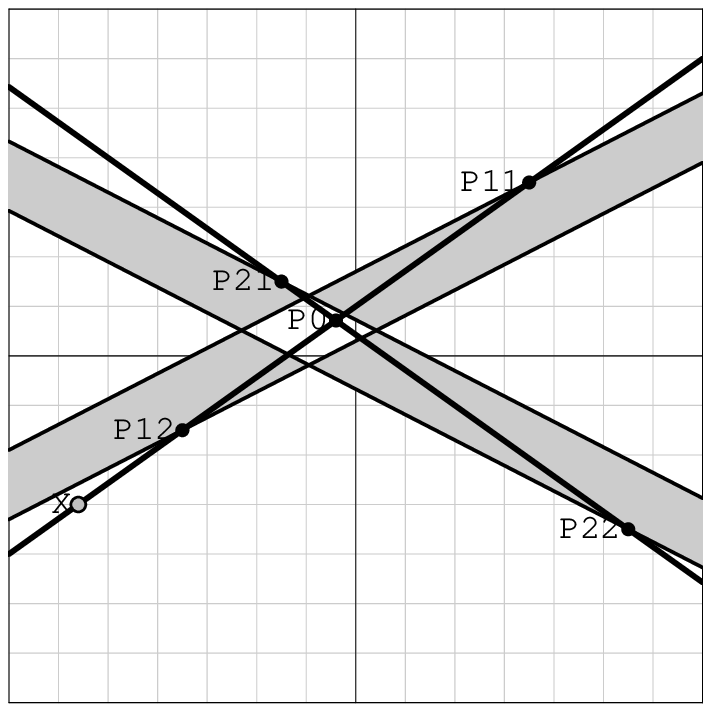}}
 \caption{} \label{c1p2}
\end{minipage}

\rule{0pt}{5pt}

\noindent
\begin{minipage}{0.49\linewidth}

\psfrag{P11}[][]{\begin{picture}(0,0)
            \put(-7,-6){\makebox(0,0)[l]{$P_{11}$}}
                        \end{picture}}
\psfrag{P12}[][]{\begin{picture}(0,0)
            \put(-4,8){\makebox(0,0)[l]{$P_{12}$}}
                        \end{picture}}
\psfrag{P21}[][]{\begin{picture}(0,0)
            \put(0,-8){\makebox(0,0)[l]{$P_{21}$}}
                        \end{picture}}
\psfrag{P22}[][]{\begin{picture}(0,0)
            \put(7,8){\makebox(0,0)[l]{$P_{22}$}}
                        \end{picture}}

\psfrag{P0}[][]{\begin{picture}(0,0)
            \put(3,8){\makebox(0,0)[l]{$P_0$}}
                        \end{picture}}

\psfrag{X}[][]{\begin{picture}(0,0)
            \put(-6,-14){\makebox(0,0)[l]{$X$}}
                        \end{picture}}

\psfrag{Q}[][]{\begin{picture}(0,0)
            \put(-3,2){\makebox(0,0)[l]{$Q$}}
                        \end{picture}}

\resizebox{\linewidth}{!}{%
  \includegraphics{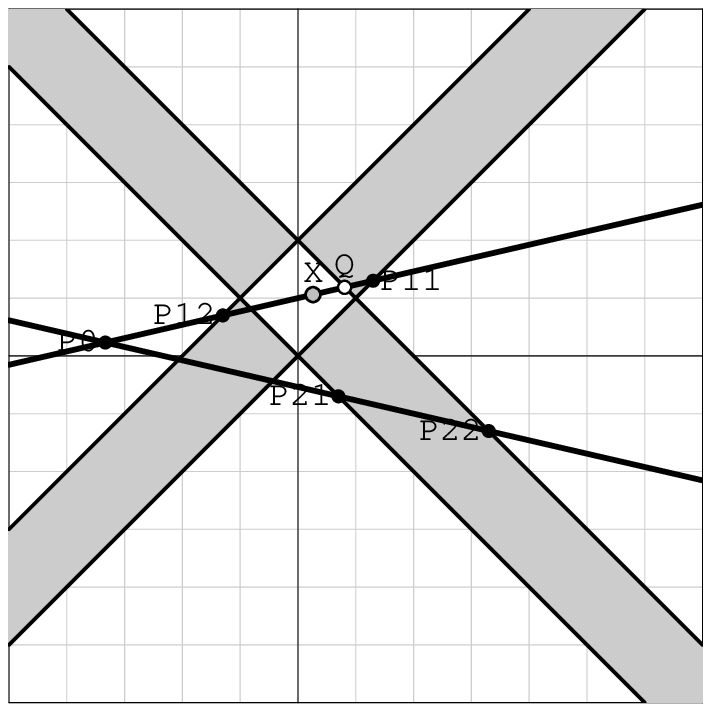}}
 \caption{} \label{c2p2}
\end{minipage}
\begin{minipage}{0.49\linewidth}

\psfrag{P11}[][]{\begin{picture}(0,0)
            \put(-6,5){\makebox(0,0)[l]{$P_{11}$}}
                        \end{picture}}
\psfrag{P12}[][]{\begin{picture}(0,0)
            \put(-5,7){\makebox(0,0)[l]{$P_{12}$}}
                        \end{picture}}
\psfrag{P21}[][]{\begin{picture}(0,0)
            \put(-3,-8){\makebox(0,0)[l]{$P_{21}$}}
                        \end{picture}}
\psfrag{P22}[][]{\begin{picture}(0,0)
            \put(13,4){\makebox(0,0)[l]{$P_{22}$}}
                        \end{picture}}

\psfrag{P0}[][]{\begin{picture}(0,0)
            \put(11,2){\makebox(0,0)[l]{$P_0$}}
                        \end{picture}}

\psfrag{X}[][]{\begin{picture}(0,0)
            \put(-8,0){\makebox(0,0)[l]{$X$}}
                        \end{picture}}

\psfrag{Q}[][]{\begin{picture}(0,0)
            \put(3,-9){\makebox(0,0)[l]{$Q$}}
                        \end{picture}}

\resizebox{\linewidth}{!}{%
  \includegraphics{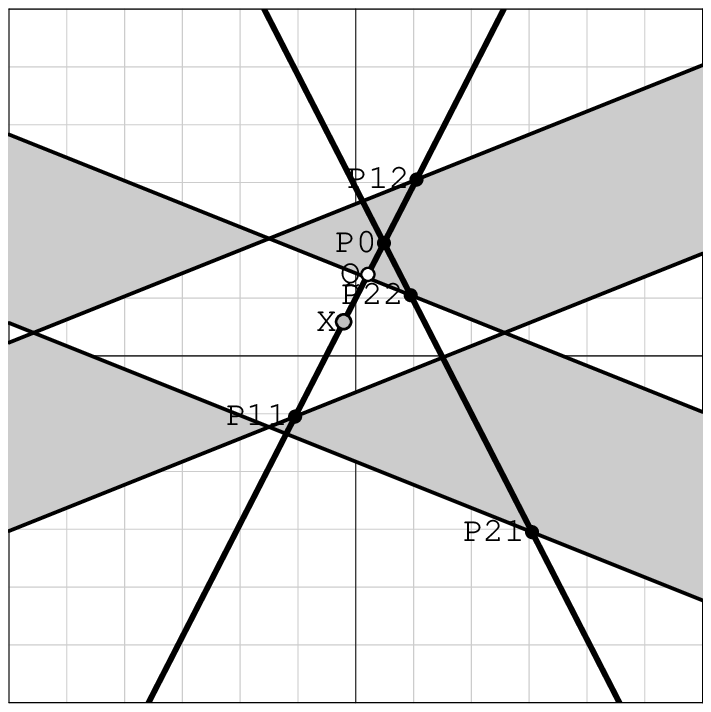}}
 \caption{} \label{c3p2}
\end{minipage}

\rule{0pt}{5pt}

\noindent
\begin{minipage}{0.49\linewidth}

\psfrag{P11}[][]{\begin{picture}(0,0)
            \put(-8,0){\makebox(0,0)[l]{$P_{11}$}}
                        \end{picture}}
\psfrag{P12}[][]{\begin{picture}(0,0)
            \put(-4,8){\makebox(0,0)[l]{$P_{12}$}}
                        \end{picture}}
\psfrag{P21}[][]{\begin{picture}(0,0)
            \put(14,0){\makebox(0,0)[l]{$P_{21}$}}
                        \end{picture}}
\psfrag{P22}[][]{\begin{picture}(0,0)
            \put(13,4){\makebox(0,0)[l]{$P_{22}$}}
                        \end{picture}}

\psfrag{P0}[][]{\begin{picture}(0,0)
            \put(12,0){\makebox(0,0)[l]{$P_0$}}
                        \end{picture}}

\psfrag{X}[][]{\begin{picture}(0,0)
            \put(-18,-5){\makebox(0,0)[l]{$X$}}
                        \end{picture}}

\psfrag{Q}[][]{\begin{picture}(0,0)
            \put(-20,4){\makebox(0,0)[l]{$Q$}}
                        \end{picture}}

\resizebox{\linewidth}{!}{%
  \includegraphics{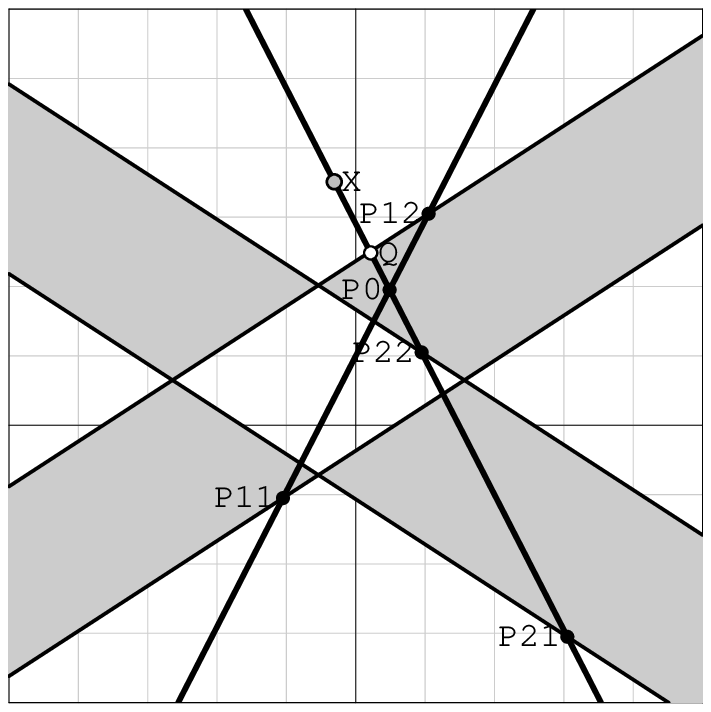}}
  \caption{} \label{c3p4}
\end{minipage}
\begin{minipage}{0.49\linewidth}

\psfrag{P11}[][]{\begin{picture}(0,0)
            \put(15,2){\makebox(0,0)[l]{$P_{11}$}}
                        \end{picture}}
\psfrag{P12}[][]{\begin{picture}(0,0)
            \put(15,2){\makebox(0,0)[l]{$P_{12}$}}
                        \end{picture}}
\psfrag{P21}[][]{\begin{picture}(0,0)
            \put(15,-3){\makebox(0,0)[l]{$P_{21}$}}
                        \end{picture}}
\psfrag{P22}[][]{\begin{picture}(0,0)
            \put(16,-3){\makebox(0,0)[l]{$P_{22}$}}
                        \end{picture}}

\psfrag{P0}[][]{\begin{picture}(0,0)
            \put(-8,0){\makebox(0,0)[l]{$P_0$}}
                        \end{picture}}

\psfrag{X}[][]{\begin{picture}(0,0)
            \put(-8,-6){\makebox(0,0)[l]{$X$}}
                        \end{picture}}

\psfrag{Q}[][]{\begin{picture}(0,0)
            \put(-7,-5){\makebox(0,0)[l]{$Q$}}
                        \end{picture}}

\resizebox{\linewidth}{!}{%
  \includegraphics{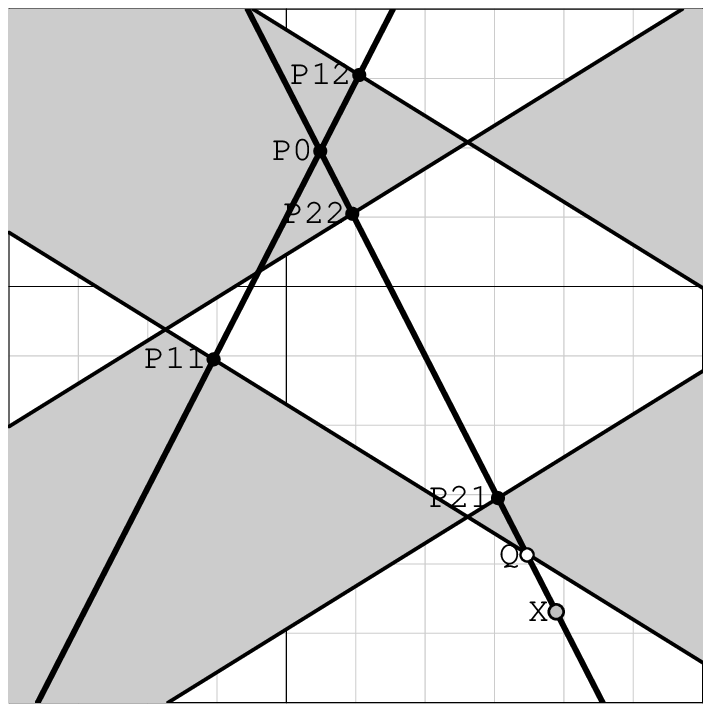}}
 \caption{} \label{c3p5}
\end{minipage}

\end{figure}

\smallskip

\noindent
{\bf Case}\!{\bf~\ref{c2}.} Points $X$ that lie in  $(\ell_1\setminus \seg_1)\cup (\ell_2\setminus \seg_2)$ are eliminated by a cross that has boundary lines through $P_{11}$ and $P_{12}$ of slope $-m$, and boundary lines through points $P_{21}$ and $P_{22}$ of slope $m$; this is similar to the method used in Figure~\ref{c1p1}.

Thus, it remains to eliminate the points inside either of the segments $\seg_1$ or $\seg_2$. Consider $X \in \seg_1$. Select a point $Q \in \ell_1$ and between $X$ and $P_{11}$ so that the line through $Q$ and $P_{22}$ has a slope $m_0 \in \nR \setminus \{-1,0\}$. The cross with boundary lines through $P_{11}$ and $P_{12}$ of slope $-m_0$, and boundary lines through points $P_{21}$ and $P_{22}$ of slope $m_0$ contains $X$ in its nucleus, thus eliminating it; see Figure~\ref{c2p2}. A similar argument works for $X\in \seg_2$.

\smallskip

\noindent
{\bf Case}\!{\bf~\ref{c3}.} We will show that no other point on $\ell_1$ and $\ell_2$, besides $P_0, P_{11}, P_{12}, P_{21},P_{22},$ belongs to all crosses in $\cC$.

Let us assume that $P_0\in \seg_1\setminus \seg_2$. Clearly, the cross that has the boundary lines through $P_{11}$ and $P_{12}$ parallel to $\ell_2$, and the boundary lines through points $P_{21}$ and $P_{22}$  parallel to $\ell_1$ eliminates any point $X$ in $\ell_1 \setminus \seg_1$, as well as the points in the interior of $\seg_2$.

Let now $X \in \seg_1$. Moreover, assume that $X$ is between $P_{11}$ and $P_0$. Select a point $Q$ on $\ell_1$ and between $X$ and $P_0$ such that the slope of the line through $Q$ and $P_{22}$ has a slope $m_0$ in $\nR \setminus \{-1,0\}$. The cross with boundary lines through $P_{11}$ and $P_{12}$ of slope $-m_0$, and boundary lines through points $P_{21}$ and $P_{22}$ of slope $m_0$ contains $X$ in its nucleus, thus eliminating it. This way, we eliminate all the points inside the segment $\seg_1$ with the exception of $P_0$; see Figure~\ref{c3p2}.

Next let $X \in \ell_2\setminus \seg_2$. First assume that $X$ and $P_0$ are on the same side of $\seg_2$. If $P_0$ is between $X$ and $\seg_2$, let $Q$ be a point on $\ell_2$ between $X$ and $P_0$ such that the slope of the line through $P_{12}$ and $Q$ is in $\nR \setminus \{0,1\}$. Denote this slope by $m_0$. The cross with boundary lines through $P_{11}$ and $P_{12}$ of slope $m_0$, and boundary lines through points $P_{21}$ and $P_{22}$ of slope $-m_0$ eliminates $X$; see Figure~\ref{c3p4}. If $X$ is between $P_0$ and $\seg_2$ a similar argument eliminates it.

If $X$ and $P_0$ are on the opposite sides of $\seg_{2}$,
then pick $Q$ on $\ell_2$ between $X$ and  $\seg_2$ such that the slope $m_0$ of the line through $P_{11}$ and $Q$ is in $\nR \setminus \{0,1\}$. The cross with boundary lines through $P_{11}$ and $P_{12}$ of slope $m_0$, and boundary lines through points $P_{21}$ and $P_{22}$ of slope $-m_0$ eliminates $X$; see Figure~\ref{c3p5}.

Finally, the case $P_0\in \seg_2\setminus \seg_1$ is handled similarly. The only exception is the case when $P_0$ is the midpoint of $\seg_2$. In this case $\tau = 1/2$ and $m = -1$ and the points in  $\seg_1\cup (\ell_2 \setminus\seg_2)$ are eliminated by crosses that have the boundary lines through $P_{11}$ and $P_{12}$ with slopes $1\pm\epsilon$, and the boundary lines through points $P_{21}$ and $P_{22}$ with slopes $-1\mp\epsilon$ for sufficiently small $\epsilon >0$.

\section{The set $\nA$}\label{secA}

In Section~\ref{secB} we proved that $\nA \subseteq \nB$ and that the set $\nB$ is a subset of the union of the sets defined in \eqref{eqsetD},  \eqref{eqsetE}, \eqref{eqsetF} and \eqref{eq6pts}.
In this section, in Subsection~\ref{srot} we prove that the set in \eqref{eqsetD} is a subset of $\nA$, in Subsection~\ref{sdeg} we prove that the set in \eqref{eqsetF} is a subset of $\nA$ and that the points in \eqref{eq6pts} are not in $\nA$, and in Subsection~\ref{srealevs} we prove that the set in \eqref{eqsetE} is a subset of $\nA$. All together these inclusions imply that $\nA$ equals the union of the sets in \eqref{eqsetD}, \eqref{eqsetE} and \eqref{eqsetF}.

\subsection{From the tetrahedron to the cone} \label{scone}

In this subsection, we show that a triple $(a,b,c) \in \cT$ if and only if $\bigl(a^2,b^2,c^2\bigr) \in \cQ$, where $\cQ$ is the interior in the first octant of the cone
\[
x^2+y^2+z^2 - 2(xy+yz+zx) = 0.
\]
That is,
\[
\cQ = \left\{ \left[\!\!
\begin{array}{c}
x \\ y \\ z
\end{array}\! \!
\right] \, : \, x,y,z > 0, \   x^2+y^2+z^2 < 2(xy+yz+zx) \right\}.
\]

By definition, $(a,b,c) \in \cT$ if and only if
\begin{equation*}
|a-b| < c < a+b.
\end{equation*}
The last two inequalities are equivalent to
\begin{equation*}
a^2 + b^2 -2ab < c^2 < a^2+b^2 + 2ab,
\end{equation*}
and these two inequalities are, in turn, equivalent to
\[
\bigl|a^2 + b^2 - c^2 \bigr| < 2 a b.
\]
Squaring both sides, followed by simple algebraic transformations,
yields the following equivalent inequalities:
\begin{align*}
\bigl( a^2 + b^2 - c^2 \bigr)^2 & < 4 a^2 b^2, \\
 a^4 + b^4 + c^4 + 2a^2b^2 - 2a^2c^2 - 2b^2c^2 & < 4 a^2 b^2, \\
 2 \bigl( a^4 + b^4 + c^4 \bigr) &
  < a^4 + b^4 + c^4 + 2 a^2 b^2 + 2a^2c^2 + 2b^2c^2,  \\
 2 \bigl( a^4 + b^4 + c^4 \bigr) &
  <  \bigl( a^2 + b^2 + c^2 \bigr)^2.
\end{align*}
Now, taking the square root of both sides, we adjust the last
inequality to look like an inequality for a dot product of two unit
vectors:
\[
\frac{a^2\cdot 1 + b^2 \cdot 1 + c^2 \cdot 1}%
{\sqrt{a^4 + b^4 + c^4} \, \sqrt{3} }
 > \sqrt{\frac{2}{3}}.
\]
This means that the cosine of the angle between the vectors $\langle a^2,b^2,c^2 \rangle$ and $\langle 1,1,1 \rangle$ is bigger than $\sqrt{2/3}$. In other words, $(a,b,c)$ are lengths of sides of a triangle if and only if the point $\bigl( a^2,b^2,c^2 \bigr)$ is inside the cone centered around the diagonal $x=y=z$ and with the angle $\arccos \sqrt{2/3}$. The algebraic equation of this cone is
\[
x^2+y^2+z^2 - 2(xy+yz+zx) = 0.
\]

The relevance of the above observation is that now we can characterize $\nA$ in terms of a ``cone preserving property'' of a class of $3\times 3$ matrices. More precisely, this observation, in combination with Subsection~\ref{1}, yields the following sequence of equivalences: \begin{align*}
(\rho,\sigma,\tau) \in \nA
 & \Leftrightarrow
    \bigl(AA_\rho,BB_\sigma,CC_\tau \bigr) \in \cT \quad \text{for all} \quad (a,b,c) \in \cT \\
 & \Leftrightarrow
    \bigl((AA_\rho)^2,(BB_\sigma)^2,(CC_\tau)^2 \bigr) \in \cQ \quad \text{for all} \quad (a^2,b^2,c^2) \in \cQ \\
  & \Leftrightarrow
  M(\rho,\sigma,\tau) \cQ \subseteq \cQ,
\end{align*}
where
\[
M(\rho,\sigma,\tau) = \left[
\begin{array}{ccc}
  \rho (\rho - 1)  & \rho  &  1-\rho  \\
  1-\sigma  &  \sigma (\sigma -1)  & \sigma  \\
 \tau  & 1-\tau  & \tau (\tau -1)
\end{array}
\right].
\]

\subsection{Rotation} \label{srot}

Consider the matrix
\[
M(\tau,\tau,\tau) = \left[
\begin{array}{ccc}
 (\tau -1) \tau  & \tau  &
   1-\tau  \\
 1-\tau  & (\tau -1) \tau  &
   \tau  \\
 \tau  & 1-\tau  & (\tau -1)
   \tau
\end{array}
\right].
\]
For simplicity, we write $M$ instead of $M(\tau,\tau,\tau)$.
A long but straightforward calculation shows that $M^T M$, that is,
\[
\left[
\begin{array}{ccc}
 (\tau -1) \tau  & 1-\tau  &
   \tau  \\
 \tau  & (\tau -1) \tau  &
   1-\tau  \\
 1-\tau  & \tau  & (\tau -1)
   \tau
\end{array}
\right] \left[
\begin{array}{ccc}
 (\tau -1) \tau  & \tau  &
   1-\tau  \\
 1-\tau  & (\tau -1) \tau  &
   \tau  \\
 \tau  & 1-\tau  & (\tau -1)
   \tau
\end{array}
\right]
\]
evaluates to $\bigl(\tau (\tau-1)+1\bigr)^2 I$ and
\[
\det M = \bigl(\tau (\tau-1)+1\bigr)^3.
\]
Therefore the matrix
\[
\bigl(\tau (\tau-1)+1\bigr)^{-1} M
\]
is orthogonal and its determinant is $1$. Thus, $\bigl(\tau (\tau-1)+1\bigr)^{-1} M$ is a rotation. Clearly, the eigenvector corresponding to the real eigenvalue is the vector $\la 1,1,1 \ra$. Therefore the transformation of $\nR^3$ induced by $M$ leaves the cone $\cQ$ invariant. Consequently, the set in \eqref{eqsetD} is a subset of $\nA$.

\subsection{Degenerate case} \label{sdeg}
Let $\tau \in \nR \setminus \{0,1\}$ and consider the matrix
\[
M = M\Bigl(\frac{1}{1-\tau}, 1-\frac{1}{\tau},\tau \Bigr) = \left[
\begin{array}{ccc}
 \dfrac{\tau }{(\tau -1)^2} & \dfrac{1}{1-\tau } & \dfrac{\tau }{\tau -1} \\[10pt]
 \dfrac{1}{\tau } & \dfrac{1-\tau }{\tau ^2} & 1-\dfrac{1}{\tau } \\[8pt]
 \tau  & 1-\tau  & (\tau -1) \tau
\end{array}
\right].
\]
A simple verification yields that the linearly independent vectors $\la 1 - \tau, 0, 1 \ra$ and $\la \tau - 1, \tau, 0 \ra$ are eigenvectors corresponding  to the eigenvalue $0$. Since $M$ has rank $1$, the third eigenvector is
\begin{equation*}
\left\la \dfrac{\tau }{\tau -1}, 1 - \dfrac{1}{\tau }, \tau (\tau-1)  \right\ra
\end{equation*}
and it corresponds to the eigenvalue
\[
\frac{\tau ^6-3 \tau ^5+3 \tau^4-\tau ^3+3 \tau ^2-3 \tau +1}%
{(\tau -1)^2 \tau ^2}.
\]
The matrix $M$ maps $\cQ$ into $\cQ$ if and only if the last eigenvector or its opposite is in $\cQ$ and the corresponding eigenvalue is positive. To explore for which $\tau$ this is the case, we introduce the change of variables
 $
\tau \mapsto (1+x)/2.
 $
Then, the eigenvector becomes
\begin{equation*}
\left\la \dfrac{x+1}{x-1}, \dfrac{x-1}{x+1}, \frac{x^2-1}{4}  \right\ra
\end{equation*}
and the corresponding eigenvalue is
\[
\frac{x^6-3 x^4+51 x^2+15}{4\left(x^2-1\right)^2}.
\]
The polynomial in the numerator is even and, since its second derivative $30 x^4 - 36 x^2+102$ is strictly positive, it is concave up. Thus, its minimum is $15$. Therefore the third eigenvalue is always positive. Next, we use the dot product to calculate the absolute value of the cosine of the angle between
\begin{equation*}
\left\la \dfrac{x+1}{x-1}, \dfrac{x-1}{x+1}, \frac{x^2-1}{4}  \right\ra
 \quad \text{and} \quad \bigl\la 1,1,1 \bigr\ra.
\end{equation*}
After simplifying, the absolute value of the cosine is
\begin{equation} \label{eqabcos}
\sqrt{ \frac{\left(x^2+3\right)^3}{3
   \bigl(x^6-7 x^4+59 x^2+11\bigr)} } .
\end{equation}
This simplification is based on the following two identities
\begin{align*}
4(x+1)^2 + 4(x-1)^2 + (x^2-1)^2 & = (3+x^2)^2 \\
4^2 (x+1)^4 + 4^2(x-1)^4 + (x^2-1)^4
 & = (3+x^2)\bigl(x^6-7 x^4+59 x^2+ 11 \bigr).
\end{align*}

Next, we need to find those values of $x$ for which~\eqref{eqabcos} is greater than $\sqrt{2/3}$. This is equivalent to
\[
\frac{\left(x^2+3\right)^3}{ x^6-7 x^4+59 x^2+11 } > 2,
\]
and, after further simplification, to
\begin{equation} \label{eqipx}
 -x^6 + 23 x^4 - 91 x^2 + 5 > 0
\end{equation}
As $5$ is a root of the polynomial $-y^3 + 23 y^2 - 91 y + 5$, factoring the last displayed polynomial yields
\begin{align*}
-x^6 + 23 x^4 - 91 x^2 + 5 & = \bigl( 5 - x^2 \bigr) \bigl( 1 - 18 x^2 + x^4 \bigr) \\
 & = \bigl( 5 - x^2 \bigr) \left( x^2 - \bigl( 9 + 4 \sqrt{5} \bigr) \right)
 \left( x^2 - \frac{1}{9 + 4 \sqrt{5}} \right)
\end{align*}
Since
\[
\sqrt{9 + 4 \sqrt{5}} = 2 + \sqrt{5} \quad \text{and} \quad \frac{1}{\sqrt{9 + 4 \sqrt{5}}} = -2 + \sqrt{5},
\]
the roots of the polynomial in~\eqref{eqipx} in the increasing order are
\[
-2-\sqrt{5}, \ \ -\sqrt{5}, \ \ 2-\sqrt{5}, \ \ -2+\sqrt{5}, \ \ \sqrt{5}, \ \ 2+\sqrt{5}.
\]
Thus, the solutions of~\eqref{eqipx} are the open intervals
\[
\left(-2-\sqrt{5}, \ -\sqrt{5} \right), \ \
 \left(2-\sqrt{5}, \ -2+\sqrt{5}\right), \ \ \left(\sqrt{5}, \ 2+\sqrt{5}\right).
\]
The corresponding intervals for $\tau$ are
\[
\left( -\phi, \ -\phi^{-1}\right), \qquad
 \left( \phi^{-2}, \ \phi^{-1}\right), \qquad
 \left( \phi,\ \phi^2 \right),
\]
with the approximate values being
\[
(-1.618, -0.618), \qquad (0.382, 0.618), \qquad (1.618, 2.618).
\]
Hence, the matrix $M$ maps $\cQ$ into $\cQ$ if and only if
\[
\tau \in \left( -\phi, \ -\phi^{-1}\right) \cup
 \left( \phi^{-2}, \ \phi^{-1}\right) \cup
 \left( \phi,\ \phi^2 \right).
\]
This, in particular implies that the set in \eqref{eqsetF} is a subset of $\nA$.

Moreover, for $\tau \in \bigl\{  -\phi, -\phi^{-1}, \phi^{-2}, \phi^{-1}, \phi,\phi^2\bigr\}$, the eigenvector of $M$ lies on the boundary of the cone $\cQ$. This boundary corresponds to the squares of the sides of degenerate triangles. Thus, for an arbitrary triangle $ABC$, for the six triples listed in \eqref{eq6pts} the corresponding cevians form a degenerate triangle. In other words, the points in \eqref{eq6pts} are not in $\nA$.

\subsection{Real eigenvalues} \label{srealevs}

Consider the matrix
\[
M = M(\xi,-\xi,2-\xi) = \left[
\begin{array}{ccc}
 (\xi -1) \xi  & \xi  & 1-\xi  \\
 \xi +1 & (\xi +1) \xi  & -\xi
   \\
 2-\xi  & \xi -1 & (1-\xi )
   (2-\xi )
\end{array}
\right].
\]
The eigenvectors of this matrix are
\[
\left[\!\!
\begin{array}{c}
-1  \\ 1 \\ 1
\end{array}\!\!
\right],
\qquad
\left[\!\!
\begin{array}{c}
1  \\ \phi^{-2} \\ \spn
\end{array}\!\!
\right], \qquad
\left[\!\!
\begin{array}{c}
1  \\ \spn \\ \phi^{-2}
\end{array}\!\!
\right],
 \quad \text{where} \quad \phi = \frac{1+\sqrt{5}}{2}.
\]
The corresponding eigenvalues are
\begin{equation*} 
-1-\xi+\xi^2, \qquad \bigl(\spn-1-\xi\bigr)^2, \qquad
  \bigl(\phi^{-2}-1-\xi\bigr)^2.
\end{equation*}
This can be verified by direct calculations. In the verification of these claims and in the calculations below, the following identities involving $\spn$ are used:
\[
(\phi^2-1)^2 = \phi^2, \quad \spn + \phi^{-2} = 3, \quad  \left(\spn - \phi^{-2}\right)^2 = 5.
\]

The matrix $M$ is singular if and only if $-1-\xi+\xi^2 = 0$. That is, for
\[
\xi = -\phi^{-1} = \phi^{-2} - 1 \quad  \text{or} \quad  \xi = \phi = \spn -1.
\]
But for $\xi = -\phi^{-1}$ the matrix $M$ of this section coincides with the matrix $M$ in Section~\ref{sdeg} with $\tau = \phi^2$ there. For $\xi = \phi$ the matrix $M$ of this section coincides with the matrix $M$ in Section~\ref{sdeg} with $\tau = \phi^{-2}$ there. Therefore in the rest of this section we can assume that $M$ is invertible, that is, we assume
\begin{equation} \label{eqassu}
\xi \neq - \phi^{-1} \quad  \text{and} \quad  \xi \neq \phi.
\end{equation}

The matrix
\[
B = \left[
\begin{array}{ccc}
-1 & 1 & 1  \\
1 &  \phi^{-2}  & \spn \\
1 &  \spn &  \phi^{-2}
\end{array}
\right]
\]
diagonalizes $M$ since
\[
B^{-1} M B = D = \left[
\begin{array}{ccc}
-1-\xi+\xi^2 & 0 & 0  \\
0 &    \bigl(\spn - 1 - \xi\bigr)^2   & 0 \\
0 &  0 &  \bigl(\phi^{-2} - 1 - \xi\bigr)^2
\end{array}
\right].
\]
Notice that only the top left diagonal entry in $D$ might be non-positive. Also, the square of the top left diagonal entry in $D$ is the product of the remaining two diagonal entries:
\begin{equation} \label{eqprevs}
\bigl(\phi^2-1-\xi\bigr)^2  \bigl(\phi^{-2}-1-\xi\bigr)^2 = \left(-1-\xi+\xi^2\right)^2.
\end{equation}

Next, we define another cone,
\[
\cQ_0 = \left\{ \left[\!\!
\begin{array}{c}
u  \\ v \\ w
\end{array}\!\!
\right] \, : \, v, w > 0, \ \  u^2 < 4 v w \right\},
\]
and prove that $D\cQ_0 = \cQ_0$ and $B \cQ_0 = \cQ$.

First notice that the definition of $\cQ_0$, \eqref{eqprevs} and~\eqref{eqassu} yield the following equivalences:
{\allowdisplaybreaks \begin{align*}
\left[\!\!
\begin{array}{c}
u  \\ v \\ w
\end{array}\!\!
\right] \in \cQ_0 \ \ & \Leftrightarrow \ \
 v, w > 0 \ \text{and} \  u^2 < 4 v w \\
 & \Leftrightarrow \ \
 \bigl(\spn - 1 - \xi\bigr)^2 v > 0, \ \ \bigl(\phi^{-2} - 1 - \xi\bigr)^2w > 0 \\
 & \rule{1cm}{0pt} \rule{0pt}{12pt}  \text{and} \ \ \\
 &  \rule{7mm}{0pt}
    \bigl(-1-\xi +\xi^2\bigr)^2 u^2\! <\! 4 \Bigl(\! \bigl(\spn - 1 - \xi\bigr)^2 v\! \Bigr)\! \Bigl(\! \bigl(\phi^{-2} - 1 - \xi\bigr)^2w \!\Bigr) \\
 & \Leftrightarrow \ \
  \left[\!\!
\begin{array}{c}
\bigl(-1-\xi +\xi^2\bigr) u  \\ \bigl(\spn - 1 - \xi\bigr)^2 v \\ \bigl(\phi^{-2} - 1 - \xi\bigr)^2 w
\end{array}\!\!
\right] \in \cQ_0 \\
  & \Leftrightarrow \ \
  D \left[\!\!
\begin{array}{c}
u  \\ v \\ w
\end{array}\!\!
\right] \in \cQ_0,
\end{align*}}
which, in turn, prove $D\cQ_0 = \cQ_0$.

Second, we verify that the matrix $B$ maps $\cQ_0$ onto $\cQ$. Set
\begin{align*}
\left[\!\!
\begin{array}{c}
x  \\ y \\ z
\end{array}\!\!
\right] = B \left[\!\!
\begin{array}{c}
u \\ v \\ w
\end{array}\!\!
\right] & = \left[\!\!
\begin{array}{c}
- u + v + w  \\
 u +  \phi^{-2} v + \spn w  \\
 u + \spn v + \phi^{-2} w
\end{array}\!\!
\right],
\end{align*}
and calculate
{\allowdisplaybreaks\begin{align*}
 x^2 + y^2 + & z^2 - 2(xy +yz +zx) \\
 & = \bigl(x^2- 2xy - 2zx \bigr) + \bigl(y^2+z^2 - 2yz  \bigr) \\
 & = x \bigl(x - 2 (y + z) \bigr) + \bigl(y-z\bigr)^2 \\
 & = (- u + v + w) \\
 & \rule{1cm}{0pt} \times
  \Bigl( - 5u + \bigl(1 - 2 ( \spn + \phi^{-2} ) \bigr)(v+w) \Bigr) \\
  & \rule{2cm}{0pt} + \bigl(\spn-\phi^{-2}\bigr)^2 (v-w)^2 \\
  & = (-5) \bigl( - u + v + w \bigr) \bigl( u + v + w \bigr) + 5 (v-w)^2 \\
   & = 5 \bigl( u^2 - (v + w)^2 + (v-w)^2 \bigr) \\
   & = 5 \bigl( u^2 - 4 v w \bigr).
\end{align*}}

Since $x^2 + y^2 + z^2 - 2(xy +yz +zx) < 0$ represents the interior (or the disconnected part) of the circular cone with the axis $x=y=z$ and with a generatrix  $x=y, z =0$, all three coordinates $x,y,z$ satisfying the last inequality must have the same sign. Clearly, if $u^2 < 4vw$, then $v,w$ are nonzero and have the same sign.

Assume
\[
x^2 + y^2 + z^2 - 2(xy +yz +zx) = 5 \bigl( u^2 - 4 v w \bigr) < 0.
\]
Then, if both $v,w < 0$, we have
\[
x = -u+v+w < 2\sqrt{|v|}\sqrt{|w|} - |v|-|w| = - \bigl(\sqrt{|v|} - \sqrt{|w|} \bigr)^2 \leq 0.
\]
The contrapositive of this implication is that, if $x > 0$, then at least one, and hence both, $v, w > 0$.

Conversely, if $v, w > 0$, then
\[
x = -u+v+w > - 2\sqrt{v}\sqrt{w} +v+w =\bigl(\sqrt{v} - \sqrt{w} \bigr)^2 \geq 0.
\]
That is: $v, w > 0$ implies $x > 0$. Hence, the following equivalences hold:
{\allowdisplaybreaks\begin{align*}
\left[\!\!
\begin{array}{c}
x  \\ y \\ z
\end{array}\!\!
\right] \in \cQ \ \  & \Leftrightarrow \ \
 x > 0 \ \text{and} \ x^2 + y^2 + z^2 - 2(xy +yz +zx) < 0 \\
 &  \Leftrightarrow \ \ v,w > 0 \ \text{and} \ u^2 < 4 vw \\
 & \Leftrightarrow \ \ \left[\!\!
\begin{array}{c}
u \\ v \\ w
\end{array}\!\!
\right] \in \cQ_0.
\end{align*}}
This proves that $B \cQ_0 = \cQ$. Since $B^{-1} M B = D$ and $D\cQ_0 = \cQ_0$, we have
\[
M \cQ = M B \cQ_0 = B D \cQ_0 = B \cQ_0 = \cQ.
\]
This proves that $M\bigl(P_{21}(2-\xi)\bigr), \xi \in \nR\setminus\{-\phi^{-1},\phi\}$, leaves the cone $\cQ$ invariant. Thus,  the second points listed in \eqref{eqsetE} belong to $\nA$. The other families of points in \eqref{eqsetE} are treated similarly, proving that the set in \eqref{eqsetE} is a subset of $\nA$.  This completes our proof of the Theorem.

\section{Closing comments}

We wish to end with a remark concerning an arbitrary {\it fixed} triangle $ABC$ with the sides $a,b,c$. For such a triangle, denote by $\nA(a,b,c)$ the set of all triples $(\rho, \sigma, \tau) \in \nR^3$ for which the cevians $AA_\rho$, $BB_\sigma$, and $CC_\tau$ form a triangle.  Clearly, $\nA \subset \nA(a,b,c)$. In fact,
 $
\nA = \bigcap_{(a,b,c) \in \mathcal T}  \nA(a,b,c) .
 $
The set $\nA(a,b,c)$ is a subset of $\nR^3$ bounded by the surface which is the union of the solutions to $AA_\rho+BB_\sigma=CC_\tau$ or $BB_\sigma+CC_\tau = AA_\rho$ or $CC_\tau+AA_\rho =BB_\sigma$. The general equation of this surface is quite cumbersome. To get an idea how this surface looks like, we plot it in Figure~\ref{SetEqA} for an equilateral triangle.
Figure~\ref{SetEqA} also illustrates the last sentence in Subsection~\ref{sdeg}: For an arbitrary triangle $ABC$, the six special points listed in \eqref{eq6pts} belong to the surface bounding
$\nA(a,b,c)$.

\section*{Acknowledgement}

The authors thank an anonymous referee for the careful reading of the paper and very useful suggestions which significantly improved the presentation.

\begin{figure}[H]
\setlength{\abovecaptionskip}{-3pt}%
\setlength{\belowcaptionskip}{-5pt}%
\psfrag{tau}[][]{\begin{picture}(0,0)
            \put(0,0){\makebox(0,0)[l]{\LARGE $\tau$}}
                        \end{picture}}
\psfrag{sigma}[][]{\begin{picture}(0,0)
            \put(0,5){\makebox(0,0)[l]{\LARGE $\sigma$}}
                        \end{picture}}
\psfrag{rho}[][]{\begin{picture}(0,0)
            \put(-5,0){\makebox(0,0)[l]{\LARGE $\rho$}}
                        \end{picture}}
\resizebox{0.8\linewidth}{!}{%
  \includegraphics{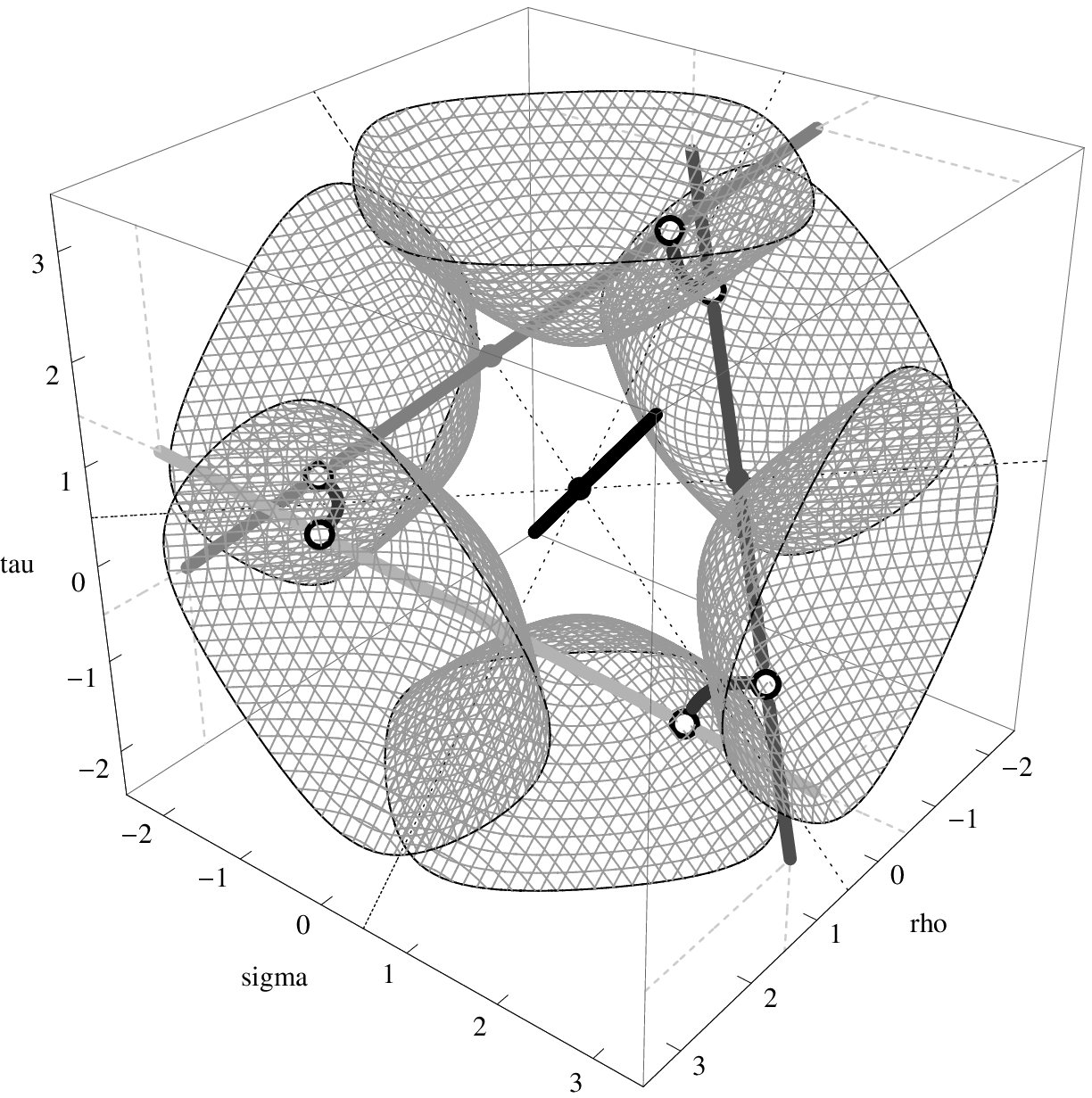}}
    \caption{The set $\nA$ and the surface bounding $\nA(1,1,1)$}
\label{SetEqA}
\end{figure}

\end{document}